\RequirePackage{fix-cm}
\documentclass[smallextended]{svjour3}       
\smartqed  
\usepackage{graphicx}
\usepackage{amssymb}
\usepackage{amsmath}
\begin{document}

\title{
On the condition numbers of a multiple generalized eigenvalue
}

\titlerunning{Condition numbers of a multiple generalized eigenvalue}        

\author{ Yuji Nakatsukasa
}


\institute{Y. Nakatsukasa \at
   Department of Mathematics, University of California, Davis, CA 95616            \\
              \email{ynakatsukasa@ucdavis.edu}           
}

\date{Received: date / Accepted: date}

\maketitle

\begin{abstract}
For standard eigenvalue problems, a closed-form expression for the condition numbers of a multiple eigenvalue is
 known. In particular, they are uniformly 1 in the Hermitian case, and generally take different values in the non-Hermitian case. 
We consider the generalized eigenvalue problem and identify the condition numbers of a multiple eigenvalue. 
Our main result is  that a  multiple eigenvalue generally has multiple condition numbers, even in the Hermitian  definite case. 
The condition numbers are characterized in terms of the singular values of the outer product of the corresponding left and right eigenvectors. 
\end{abstract}
\section{Introduction}
This paper is concerned with the sensitivities of a 
 multiple eigenvalue in generalized eigenvalue problems $Ax=\lambda Bx$. 
When perturbations are introduced in the matrices, a multiple eigenvalue $\lambda_0$ of multiplicity $r$ generally splits into $r$ simple eigenvalues. We are interested in how sensitively each of the $r$ eigenvalues can be perturbed. 

The first order perturbation expansion for a simple eigenvalue is a well-known result \cite{stewart_2,stewart-sun:1990}, and that for a multiple eigenvalue is also studied in \cite{sun92} for standard eigenvalue problems, and in \cite{1storder} for generalized eigenvalue problems, including singular matrix pairs. 
Using such results and following the definition of $r$ condition numbers of a multiple eigenvalue for standard eigenvalue problems introduced in \cite{sunnm95}, we naturally define condition numbers for the generalized case, as shown below in \eqref{defcond}. 

Sun shows in \cite{sunnm95} that the $r$ condition numbers $ \kappa_i(A,\lambda_0)$ for $i=1,\ldots, r$
of a nondefective multiple eigenvalue of a non-Hermitian matrix are expressed by
\begin{equation}  \label{sunresult}
 \kappa_i(A,\lambda_0)=\left(\prod_{j=1}^ic_j(A,\lambda_0)\right)^{1/i} \quad\mbox{for}\quad i=1,\ldots,r,
\end{equation}
where $c_j(A,\lambda_0)$  are the secants of the canonical angles between the left and right invariant subspaces corresponding to the multiple eigenvalue. When $A$ is non-Hermitian $c_i(A,\lambda_0)$ generally take different values for different $i$, hence so do $\kappa_i(A,\lambda_0)$ and \eqref{sunresult} shows that a multiple eigenvalue has multiple condition numbers. 
Contrarily, in the Hermitian case we have $c_j(A,\lambda_0)\equiv 1$, so  $\kappa_i(A,\lambda_0)=1$ for $i=1,\ldots,r$. Hence \eqref{sunresult} also shows the well-known fact that a multiple eigenvalue of a Hermitian matrix always has a uniform condition number $1$.

Since a standard  non-Hermitian eigenvalue problem can be regarded as a special case of the generalized non-Hermitian eigenvalue problem  $Ax=\lambda Bx$, it is clear that a multiple eigenvalue in this case must have multiple condition numbers. 
On the other hand, in the important case of the generalized Hermitian definite pair (where $A$ and $B$ are Hermitian and $B$ is positive definite), it is not trivial whether or not the condition numbers $\kappa_i(A,B,\lambda_0)$ for $i=1,\ldots, r$ take different values. 
In this paper we identify their expression, which shows that they generally do take $r$ different values.
We shall see that there are two sources for this different conditioning, namely the difference between the left and right eigenvectors (as present in non-Hermitian standard eigenproblems), and the fact that the $B$-orthonormal eigenvectors have different 2-norms (as present in the case $B\not= I$). 

It is important to note that in the Hermitian definite case, a natural choice of metric can be based on the $B$-based inner product $(x,y)_B=x^HBy$ instead of the standard inner product $(x,y)=x^Hy$. 
This leads to the 
standard eigenvalue problem for the Hermitian matrix $B^{-1/2}AB^{-1/2}$, so in this inner product all the condition numbers are the same (see the second remark in section \ref{seciden}). 
 Our discussion in this paper assumes the use of the standard inner product. 


 The condition numbers $\kappa_i$ for $i=1,\ldots, r$ of a nondefective finite multiple eigenvalue in four situations are summarized in Table \ref{sumresults}, expressed in terms of $\sigma_1,\ldots,\sigma_r$, the $r$ positive singular values of $X_1Y_1^H$ (see section \ref{conddef}).  The contribution of this paper is that we fill in the second row, 
that is, we  identify the condition numbers of a multiple eigenvalue in a generalized eigenvalue problem, both for the Hermitian definite case and the non-Hermitian case. 
Here $\tau$ is a prescribed positive constant that accounts for perturbation scalings, see section \ref{conddef}. 

\begin{table}[htbp]
\caption{Summary of condition numbers $\kappa_i$ of a multiple eigenvalue for $i=1,\ldots, r$.}
\label{sumresults}
\begin{center}
    \begin{tabular}{|c|c|c|}\hline
               & Hermitian & Non-Hermitian\\\hline
$Ax=\lambda x$ &1 & $\left(\prod_{j=1}^i\sigma_j\right)^{1/i}$ \\\hline
$Ax=\lambda Bx$&$(1+\tau|\lambda_0|)\min_{1\leq j\leq i}\sqrt{{\sigma_j\sigma_{i-j+1}}}$ & $(1+\tau|\lambda_0|)\left(\prod_{j=1}^i\sigma_j\right)^{1/i}$\\\hline
    \end{tabular}
\end{center}
  \end{table}
There are a number of related studies in the literature. \cite{holder09} investigates the H$\ddot{\rm o}$lder condition number, which is essentially $\kappa_1$ in our terminology when $\lambda_0$ is nondefective. The focus of \cite{holder09} is the effect of the structure of the perturbation on the H$\ddot{\rm o}$lder condition number, and  in section \ref{structure} we discuss how our results are related to those in \cite{holder09}. 

An observation that a multiple generalized eigenvalue has different sensitivities under perturbations was first made in \cite[p.300]{stewart-sun:1990}, which mentions that a multiple eigenvalue of a pair such as $A=  \begin{bmatrix} 2000&0\\0&2  \end{bmatrix}$, $B=  \begin{bmatrix} 1000&0\\0&1  \end{bmatrix}$ 
tends to behave differently under perturbations in $A$ and $B$. We note that as shown in \cite{Veselic1993}, for Hermitian definite pairs, small componentwise relative changes in $A$ and $B$ can introduce only small relative perturbation to any eigenvalue, and it is easy to see the two eigenvalues of the above pair $(A,B)$ have similar perturbation behaviors. However, in terms of ``standard'' normwise perturbation, that is, when $(A,B)$ is perturbed to $(A+\epsilon E,B+\epsilon F)$ under $\|E\|_2,\|F\|_2\leq 1$ and $\epsilon\rightarrow 0$, a multiple eigenvalue can exhibit different behaviors. 
\cite{mybit,limetruhar} consider the Hermitian definite case and give an explanation for this behavior, presenting $r$ different perturbation bounds for $\lambda_0$ under perturbations of finite norm. 
The approach of this paper is different in that our focus is on the condition numbers, which are attainable perturbation bounds in the first order sense in the limit $E,F\rightarrow 0$. The bounds in \cite{mybit,limetruhar} are valid for non-asymptotic $E,F$ but are less tight (generally not attainable) when  $E,F\rightarrow 0$. 

Our arguments closely follow that of \cite{sunnm95}, in which the condition numbers are called \emph{worst-case} condition numbers, to emphasize the difference from the \emph{typical-case}  condition numbers, as presented in \cite{stewartgraded}. 
In this sense, our results should also be regarded as worst-case condition numbers, in that $\kappa_i$  are the largest attainable bounds in the first order sense. Experiments show that these bounds are not likely to be attained in practice for randomly generated perturbations, especially for large $i$ (see the example in section \ref{structure}). 

The structure of this paper is as follows. First in section \ref{prelim} we establish a characterization of condition numbers of a multiple generalized eigenvalue. 
In section \ref{ghep} we treat the Hermitian definite case and show that a multiple generalized eigenvalue indeed has multiple condition numbers. Section \ref{gnep} treats the general non-Hermitian case, 
in which we give simple numerical examples to illustrate our results. 
We also discuss the case of an infinite and defective multiple eigenvalue.

Notations: $\sigma_i(X)$ denotes the $i$th largest singular value of a rectangular matrix $X$. $\lambda_i(A)$ denotes the $i$th eigenvalue of $A$, arranged in decreasing order of magnitude. 
We use  only the matrix spectral norm $\|\cdot\|_2$, so that $\|A\|_2=\sigma_1(A)$. We use MATLAB notation, in which $V(i,j:k)$ denotes the $j$th to $k$th elements of the $i$th row of $V$. $I_k$ is the identity matrix of order $k$. 
\section{Condition numbers of a multiple generalized eigenvalue}\label{prelim}
\subsection{Definition}\label{conddef} 
For an $n$-by-$n$ matrix pair $(A,B)$, suppose that $\lambda_0$ is a nondefective finite multiple eigenvalue (we discuss the infinite and defective cases later in section \ref{defsec}) of multiplicity $r$, so that there exist nonsingular matrices $X=(X_1,X_2)$ and $Y=(Y_1,Y_2)$ with $X_1,Y_1\in \mathbb{C}^{n\times r}$ that satisfy 
\begin{equation}  \label{xy}
Y^HAX=\begin{bmatrix} \lambda_0 I_r&0\\0&J_A \end{bmatrix},\quad Y^HBX=\begin{bmatrix}  I_r&0\\0&J_B \end{bmatrix}. 
\end{equation}
Here the spectrum of the pair $(J_A,J_B)$ does not contain $\lambda_0$.
Then, the pencil $(A+\epsilon E,B+\epsilon F)$ has eigenvalues $\widehat\lambda_1,\widehat\lambda_2, \ldots,\widehat\lambda_r$ admitting the first order expansion \cite{moro97,holder09}
\begin{equation}  \label{expand}
\widehat \lambda_i= \lambda_0 +\lambda_i(Y_1^H(E-\lambda_0F)X_1)\epsilon+o(\epsilon), \quad i=1,2,\ldots, r, 
\end{equation}
where $\lambda_i(Y_1^H(E-\lambda_0F)X_1)$ are the eigenvalues of $Y_1^H(E-\lambda_0F)X_1$ for $i=1,\ldots,r$. 
In light of \eqref{expand} and following the definition presented in \cite{sunnm95}, we define $r$ condition numbers $\kappa_i(A,B,\lambda_0)$ for $i=1,\ldots,r$ of the multiple eigenvalue $\lambda_0$ as follows. 
  \begin{definition}\label{def}
Let an $n$-by-$n$ matrix pair $(A,B)$ have the decomposition \eqref{xy}, and let $\tau>0$ be a prescribed constant. 
We define the  condition numbers of  $\lambda_0$, a multiple eigenvalue of $(A,B)$ of multiplicity $r$, by
\begin{equation}  \label{defcond}
\kappa_i(A,B,\lambda_0)\equiv \sup_{\|E\|_2\leq 1,\|F\|_2\leq \tau}|\lambda_i(Y_1^H(E-\lambda_0F)X_1)|, \quad i=1,\ldots, r,
\end{equation}
where the eigenvalues $\lambda_i(Y_1^H(E-\lambda_0F)X_1)$ are 
ordered such that $|\lambda_1(Y_1^H(E-\lambda_0F)X_1)|\geq |\lambda_2(Y_1^H(E-\lambda_0F)X_1)|\geq \cdots \geq |\lambda_r(Y_1^H(E-\lambda_0F)X_1)|.$
  \end{definition}

In words, $\kappa_i(A,B,\lambda_0)$ measures by how much small changes in $A$ and $B$ can be magnified in the multiple eigenvalue $\lambda_0$, in the first order sense. 
$\tau$ is a positive constant that allows for the case where perturbations in $A$ and $B$ occur in different magnitudes, which is a notion adopted for example in \cite{highig}.

\subsection{Equivalent characterization of $\kappa_i(A,B,\lambda_0)$}
Here we show that  $\kappa_i(A,B,\lambda_0)$ can be expressed in an equivalent form, just as in the standard ($B=I_n$) case $\kappa_i(A,\lambda_0)$ can be expressed as \eqref{sunresult} using the secants of the canonical angles $c_j(A,\lambda_0)$ between the left and right invariant subspaces corresponding to $\lambda_0$. 
Note that  
 $c_j(A,\lambda_0)=\sigma_j(X_1Y_1^H)$ for $j=1,\ldots,r$. 
In this paper we use the quantity $\sigma_j(X_1Y_1^H)$ instead of the canonical angles to identify the condition numbers, because   it allows us to treat generalized eigenvalue problems in a uniform way. 

We use the proof of Theorem 2.1 in \cite{sunnm95}, whose crucial identity is
\begin{align*}
|\lambda_i(Y_1^H(E-\lambda_0F)X_1)|&=|\lambda_i((E-\lambda_0F)X_1Y_1^H)|=|\lambda_i((E-\lambda_0F)U\Sigma V^H)|\\
&=| \lambda_i( V^H(E-\lambda_0F)U\Sigma)| 
\end{align*}
for $i=1,\ldots,r$, where  $X_1Y_1^H=U\Sigma V^H$ is the ``thin'' SVD. Here, to get the first and last equalities we used  the fact \cite[p.27]{stewart-sun:1990} that for general $X\in \mathbb{C}^{n\times m}$ and $Y\in \mathbb{C}^{m\times n}$, the nonzero eigenvalues of $XY$ and those of $YX$ are the same. 
Since $V$ and $U$ have orthonormal columns and $E,F$ can take arbitrary matrices with $\|E\|_2\leq 1,\|F\|_2\leq \tau$, 
it follows that $V^HEU, V^HFU$ can also take arbitrary matrices such that $\|V^HEU\|_2\leq 1, \|V^HFU\|_2\leq \tau$. 
Hence, redefining $E:=V^HEU$ and  $F:=V^HFU$, we see that the condition numbers $\kappa_i(A,B,\lambda_0)$ have the following  equivalent  characterization. 
\begin{lemma} Under the assumptions in Definition \ref{def}, suppose that $X_1Y_1^H=U\Sigma V^H$ is the SVD where $\Sigma=\mbox{diag}(\sigma_i)$ is $r$-by-$r$ $(\sigma_1\geq\sigma_2\geq \ldots \geq \sigma_r> 0)$. Then, $\kappa_i(A,B,\lambda_0)$ in \eqref{defcond} can be expressed as
\begin{equation}  \label{defcond2}
\kappa_i(A,B,\lambda_0)\equiv \sup_{\|E\|_2\leq 1,\|F\|_2\leq \tau}|\lambda_i(\Sigma(E-\lambda_0F))|, \quad i=1,\ldots, r.
\end{equation}  
\end{lemma}
Here we have $\sigma_r>0$ because both $X_1$ and $Y_1$ have full column-rank. 
Note that the size of $E$ and $F$ in \eqref{defcond2} is $r$-by-$r$, which is smaller than $n$-by-$n$ as in \eqref{defcond}. 
In the sequel we use the expression \eqref{defcond2} to identify $\kappa_i(A,B,\lambda_0)$. 
\section{Hermitian definite pairs}\label{ghep}
\subsection{Specifications}
When $(A,B)$ is a Hermitian definite pair, all the eigenvalues are always real and nondefective, and there exists a nonsingular matrix $X$ such that \cite{golubbook}
\begin{equation}  \label{x}
X^HAX=\begin{bmatrix} \lambda_0 I_r&0\\0&\Lambda_1 \end{bmatrix},\quad X^HBX=I_n, 
\end{equation}
where $\Lambda_1$ is a diagonal matrix containing the eigenvalues not equal to $\lambda_0$. 
Hence the diagonals of $\Sigma$ in \eqref{defcond2} are the $r$ positive singular values of the matrix $X_1X_1^H$, which are equal to the eigenvalues of the matrix $X_1^HX_1$. 
Since  $(A,B)$ is a Hermitian definite pair it is natural to require that the  perturbation matrices preserve the property, so  \eqref{defcond2} becomes the  ``structured'' condition numbers $\kappa_i(A,B,\lambda_0;\mathbb{S})$, 
 expressed by 
\begin{equation}  \label{defcond3}
\kappa_i(A,B,\lambda_0;\mathbb{S})\equiv \sup_{\|E\|_2\leq 1,\|F\|_2\leq \tau\atop E=E^H,F=F^H}|\lambda_i(\Sigma(E-\lambda_0F))|, \quad i=1,\ldots, r.
\end{equation}

Denoting $D=\Sigma^{1/2}=\mbox{diag}(\sqrt{\sigma_1},\ldots,\sqrt{\sigma_r})$, we see that the eigenvalues of  $\Sigma(E-\lambda_0F)$ are equal to those of the Hermitian matrix $D(E-\lambda_0F)D$. 

We further observe that $E-\lambda_0F$ can represent an arbitrary Hermitian matrix $H$ with $\|H\|_2\leq 1+\tau|\lambda_0|$, which can be done by letting $E=H/\|H\|_2$ and $F=-\tau E|\lambda_0|/\lambda_0$. Conversely, it is easily seen that the class of Hermitian matrices $H$ with $\|H\|_2\leq 1+\tau|\lambda_0|$ includes all the matrices expressed by $E-\lambda_0F$. 
Together with the fact that the singular values of a Hermitian matrix are simply the absolute values of the eigenvalues,  we have yet another characterization of condition numbers in the Hermitian definite case
\begin{equation}  \label{defcond4}
\kappa_i(A,B,\lambda_0;\mathbb{S})= (1+\tau|\lambda_0|)\sup_{\|H\|_2\leq 1\atop H=H^H}\sigma_i(DHD), \quad i=1,\ldots, r. 
\end{equation}
\subsection{Identifying the condition numbers}\label{seciden}
Now we are ready to identify the condition numbers $\kappa_i(A,B,\lambda_0;\mathbb{S})$ using the expression  \eqref{defcond4}. 
\begin{theorem}\label{mainthm}
In the Hermitian definite case, $\kappa_i(A,B,\lambda_0;\mathbb{S})$ as in \eqref{defcond3},\eqref{defcond4} is
\begin{equation}
  \label{main}
\kappa_i(A,B,\lambda_0;\mathbb{S})= (1+\tau|\lambda_0|)\min_{1\leq j\leq i}\sqrt{{\sigma_j\sigma_{i-j+1}}},\quad i=1,\ldots,r. 
\end{equation}
\end{theorem}
Note that 
 $\sqrt{{\sigma_j\sigma_{i-j+1}}}$ is  the geometric mean of $\sigma_j$ and $\sigma_{i-j+1}$, the $j$th largest and smallest of  $(\sigma_1,\sigma_2,\ldots, \sigma_i)$, which is the set of the  $i$ largest singular values of $X_1^HX_1$. 
\bigskip

$Proof$  
In view of \eqref{defcond4}, to prove the theorem it suffices to prove  that for any Hermitian $H$ such that $\| H\|_2=1$, 
$\sigma_i(DHD)$ 
 is bounded  above by $\min_j\sqrt{{\sigma_j\sigma_{k-j+1}}}$ for $i=1,\ldots,r$, and that this bound is attainable. 

First, proving attainability is simply done by considering the case where $H$ is zero except for its $i\times i$ leading principal submatrix, which is set to the antidiagonal matrix (which has $1$ on the antidiagonals and 0 elsewhere). This choice of  $ H$ makes the $i\times i$ leading principal submatrix of $D HD$ also an anti-diagonal matrix, whose $j$th antidiagonal is $\sqrt{\sigma_j\sigma_{i-j+1}}$. Hence we have $\sigma_i(DHD)=\min_j\sqrt{\sigma_j\sigma_{i-j+1}}$. 

Our remaining task is to  prove that  $\min_j\sqrt{{\sigma_j\sigma_{i-j+1}}}$ is an upper bound of $\sigma_i(DHD)$ for any Hermitian $H$ with $\|H\|_2\leq 1$. 
Using the max-min characterization of singular values  \cite[p. 68]{stewart1}, we have 
\[\sigma_{i}(DH D)=\max_{Q^HQ=I_i}\min_{\|v\|_2=1}\|DH DQv\|_2,\]
 so it suffices to show that for any $Q\in\mathbb{C}^{r\times i}$ with orthonormal columns, there exists a unit vector $v$ such that  $\| DH DQv\|_2\leq \min_j\sqrt{\sigma_j\sigma_{i-j+1}}$. 

To prove this, let $j_0=\mbox{argmin}_{j\leq (i+1)/2}\sqrt{{\sigma_j\sigma_{i-j+1}}}$. Since for any $Q$ we have $\mbox{rank}(Q(1:i-j_0,:))\leq i-j_0$, 
there are at least $j_0$ linearly independent vectors in $\mathbb{C}^{i\times 1}$ that are orthogonal to the rows of $Q$. 
Therefore there must exist $P\in\mathbb{C}^{i\times j_0}$ with orthonormal columns such that the first $i-j_0$ rows of the $r$-by-$j_0$ matrix $QP$ 
are all zeros. For such $P$,  we have 
\[\|DQP\|_2=\|\mbox{diag}(0,\ldots,0,\sqrt{\sigma_{i-j_0+1}},\ldots,\sqrt{\sigma_{r}})QP\|_2\leq \sqrt{\sigma_{i-j_0+1}}.\]
 Furthermore, since $\mbox{rank}(HDQP(1:j_0-1,:))\leq j_0-1$, there must exist a unit vector $w\in\mathbb{C}^{j_0\times 1}$
that is orthogonal to $HDQP(1:j_0-1,:)$, so that 
the first $j_0-1$ rows of $HDQPw$ are all zeros. We easily see that for such $w$ we have $\|DHDQPw\|_2\leq \sqrt{\sigma_{j_0}\sigma_{i-j_0+1}}$. 
Therefore we have shown that for any $Q\in\mathbb{C}^{k\times i}$ with orthonormal columns there exists a unit vector $v_0=Pw$ such that 
\begin{align*}
  \min_{\|v\|_2=1}\|DHDQv\|_2&\leq \|DHDQv_0\|_2\leq \sqrt{\sigma_{j_0}\sigma_{i-{j_0}+1}}=\min_j\sqrt{{\sigma_i\sigma_{i-j+1}}}. 
\end{align*}
\hfill$\square$

Three remarks are in order. 
\begin{itemize}
\item When $B\not= I_n$, $\sigma_i$ for $i=1,\ldots,r$ generally take different values, so \eqref{main} shows that a multiple generalized eigenvalue has multiple condition numbers, which is our main result. Note 
that the ratio among the condition numbers  is bounded by $\kappa_1(A,B,\lambda_0;\mathbb{S})/\kappa_r(A,B,\lambda_0;\mathbb{S})\leq \sigma_1/\sigma_r$. 
Now since $\sigma_{\min}(B^{-1})\leq \sigma_r\leq\sigma_1\leq \sigma_{\max}(B^{-1})$, 
we have $\sigma_1/\sigma_r\leq \sigma_{\max}(B)/\sigma_{\min}(B)=\kappa_2(B)$, the standard $2-$norm condition number of $B$. 
It follows that if $B$ is well-conditioned then a multiple eigenvalue of a Hermitian definite pair must have similar condition numbers. 
\item For standard Hermitian eigenvalue problems $(B=I_n)$, we have $d_i\equiv 1$ and $\tau=0$, so \eqref{main} reduces to $\kappa_i=1$ for all $i$, the well-known result that a multiple eigenvalue of a Hermitian matrix has a uniform condition number 1. When one allows for perturbation in $B=I_n$, the condition numbers are $1+\tau|\lambda_0|$ regardless of the multiplicity. 
The second term suggests that larger changes occur in larger eigenvalues. This observation can be summarized as follows: for perturbation in $A$, all the eigenvalues have the same sensitivity in the absolute sense, while for perturbation in $B$, all the eigenvalues have the same sensitivity in the relative sense. 
\item The above arguments show that the difference among condition numbers of a multiple eigenvalue is due to the difference among the $r$ singular values of $X_1Y_1^H$, the outer product of the left and right eigenvectors corresponding to $\lambda_0$. $\sigma_i(X_1Y_1^H)$ are all 1 in the standard Hermitian  case because $X_1=Y_1$ and it is an orthogonal matrix. In the standard non-Hermitian case $X_1\neq Y_1$ and neither is orthogonal, so $\sigma_i(X_1Y_1^H)$ take $r$ different values. In the generalized Hermitian case we have $X_1=Y_1$ but $X_1$ is not orthogonal, so $\sigma_i(X_1X_1^H)$ again take $r$ different values. Note that one uses the $B$-based inner product this difference disappears because $X_1$ is $B$-orthogonal, recall the remark in the introduction. 
\end{itemize}

\section{Non-Hermitian pairs}\label{gnep}
Here we consider the case where $(A,B)$ is a general non-Hermitian pair. 
In view of \eqref{defcond2}, our task is to bound  $|\lambda_i(X\Sigma)|$ for an arbitrary square matrix $X$ such that $\|X\|_2\leq 1$. 
This is in fact the exact same problem addressed in \cite[Thm.3.1]{sunnm95}. Hence  the analysis there can be directly applied to yield the following result. 
\begin{theorem}
  For a non-Hermitian pair $(A,B)$ that satisfies \eqref{xy}, let  $\sigma_1\geq\sigma_2\geq \ldots \geq \sigma_r> 0$  be the positive singular values of the matrix $X_1Y_1^H$. Then, $\kappa_i(A,B,\lambda_0)$ in \eqref{defcond2} can be expressed as 
\begin{equation}  \label{main2}
\kappa_i(A,B,\lambda_0)= (1+\tau|\lambda_0|)\left(\prod_{j=1}^i\sigma_j\right)^{1/i}, \quad i=1,\ldots,r.
\end{equation}
\end{theorem}
\subsection{Structured perturbation}\label{structure}
It is instructive to revisit the Hermitian definite case, but now allowing for non-Hermitian perturbations, that is,  $E,F$ are general matrices whose norms are bounded by 1. In this case, the condition numbers $\kappa_i(A,B,\lambda_0)$ have the characterization \eqref{defcond2} (instead of \eqref{defcond4}), so they have the expression \eqref{main2}, the same as that for the non-Hermitian pair. 

As might be expected, 
 the condition number under Hermitian perturbation \eqref{main} is always no larger than that under a non-Hermitian perturbation \eqref{main2}: 
\begin{align*}
\frac{\kappa_i(A,B,\lambda_0;\mathbb{S})}{1+\tau|\lambda_0|}&= \left(\min_{1\leq j\leq i}(\sigma_j\sigma_{i-j+1})^{i}\right)^{1/2i}
\leq  \left(\prod_{j=1}^i (\sigma_j\sigma_{i-j+1})\right)^{1/2i}\\
&=  \left(\prod_{j=1}^i \sigma_j^2\right)^{1/2i}
= \left(\prod_{j=1}^i\sigma_j\right)^{1/i}=\frac{\kappa_i(A,B,\lambda_0)}{1+\tau|\lambda_0|}.
\end{align*}
The above arguments imply that if the singular values of $X_1Y_1^H$ are the same, then  under a general non-Hermitian perturbation the condition numbers $\kappa_i(A,B,\lambda_0)$ are all the same\footnote{In fact the entire first order perturbation expansions become the same.}, regardless of whether $(A,B)$ is Hermitian definite or  non-Hermitian. Therefore, the structure of the perturbation matrices plays an important role in the perturbation sensitivity of a multiple generalized eigenvalue. 
We note that the standard Hermitian case with $B\equiv I$ is an exception, in which the condition numbers are always all $1$ whether or not the perturbation matrices are Hermitian. 

This point of view, to focus on the effect of the structure of the perturbation, was investigated extensively in \cite{holder09}, in which (Theorem 4.5) it is shown that (among other structures they consider) the Hermitian structure of the perturbation matrices does not have any effect on the H$\ddot{\rm o}$lder condition number. 

At first sight this seems to contradict  our results, which show that the Hermitian structure of the perturbation matrices does affect the condition numbers of the multiple eigenvalue $\lambda_0$. The explanation is that 
\cite{holder09} treats only the H$\ddot{\rm o}$lder condition number, which is 
equivalent to $\kappa_1(A,B,\lambda_0)$  in the nondefective case. 
Here we are identifying individual condition numbers of each of the  $r$ eigenvalues. In fact, we can see that  for $i=1$,  $\kappa_i$ in \eqref{main} and \eqref{main2} are the same, both equal to $(1+\tau|\lambda_0|)\sigma_1$. 
We can easily see that they are equal also for $i=2$. The difference between \eqref{main} and \eqref{main2}  starts to take effect only for $i\geq 3$, so $\lambda_0$'s multiplicity $r$ must be at least 3. 
In particular, for a simple eigenvalue the Hermitian structure of the perturbation has no effect on the condition number, which is a trivial consequence of the results in \cite{holder09}. 
\subsection{Examples}\label{ex1}
Here we present two simple examples to illustrate the above results and observations. 
\begin{example}
 For the Hermitian definite pair $A=  \begin{bmatrix} 2000&0\\0&2  \end{bmatrix}$, $B=  \begin{bmatrix} 1000&0\\0&1  \end{bmatrix}$ presented in \cite[p.300]{stewart-sun:1990}, we have $\kappa_1(A,B,\lambda_0;\mathbb{S})=\kappa_1(A,B,\lambda_0)=3$ and $ \kappa_2(A,B,\lambda_0;\mathbb{S})=\kappa_2(A,B,\lambda_0)=3/\sqrt{1000}$, which explains why the multiple eigenvalue $\lambda_0=2$ has different sensitivities. Note that in this case the structure of the perturbation has no effect on the condition numbers, because the multiplicity of $\lambda_0$ is $r< 3$. 
\end{example}
\begin{example}
We consider a $4$-by-$4$ Hermitian definite pair $(A,B)$ expressed by 
\[A=W^H\Lambda W,\quad B=W^HW,\]
where $\Lambda=\mbox{diag}(1,1,1,2)$ and $W=\mbox{diag}(1,2,100,1)$, so the eigenvalues of $(A,B)$ are $1,1,1,2$. 
Since $X$ that diagonalizes $A,B$ (as in \eqref{x}) is $X=W^{-1}=\mbox{diag}(1,0.5,0.01,1)$ and $X_1$  is its first three columns, 
the singular values of $X_1X_1^H$ are $\sigma_1=1^2,\sigma_2=0.5^2,\sigma_3=0.01^2$ (where in this example we let $\tau=1$),  
hence by \eqref{main} it follows that $\kappa_1(A,B,1;\mathbb{S})=2,\kappa_2(A,B,1;\mathbb{S})=1$ and $\kappa_3(A,B,1;\mathbb{S})=0.02$. 
Using MATLAB version 7.10 we generated  $10^6$ sets of  random Hermitian perturbation matrices $E$ and $F$ such that  $||E||_2,||F||_2\leq 1$, and examined the behavior of the three eigenvalues of the pair $(A+\epsilon E,B+\epsilon F)$ that are closest to $\lambda_0= 1$, where we let $\epsilon=10^{-5}$. 
Specifically, denoting by $\widehat \lambda_i$ for $i=1,2,3$ the three eigenvalues of $(A+\epsilon E,B+\epsilon F)$ that are closest to 1 such that $|\widehat \lambda_1-1|\geq |\widehat \lambda_2-1|\geq |\widehat \lambda_3-1|$, we examine how large $|\widehat \lambda_i-1|/\epsilon$ can be for $i=1,2,3$.  

We also experimented with non-Hermitian perturbations, in which case we let $E,F$ be arbitrary non-Hermitian matrices with $||E||_2,||F||_2\leq 1$.  In this case the condition numbers \eqref{main2} are   $\kappa_1(A,B,1)=2,\kappa_2(A,B,1)=2(1\cdot 0.5^2)^{1/2}=1$ and
 $\kappa_3(A,B,1)=2(1\cdot 0.5^2\cdot 0.01^2)^{1/3}\simeq 0.058$, in which we confirm that the first two are the same as in the above Hermitian case. 

Lastly, in order to see how the Hermitian property of the matrices plays a role in the eigenvalue perturbation behaviors, we also tested with a non-Hermitian pair $(A,B)$ that has the same eigenvalues and $\sigma_i$ (of $X_1Y_1^H$) as the above Hermitian pair. 
We formed such a pair $(A,B)$ by defining $A=Y^{-H}\Lambda X^{-1}$ and $B=Y^{-H}X^{-1}$, where $\Lambda=\mbox{diag}(1,1,1,2)$, $Y_1^H$ (the first $3$ rows of $Y$) is set to $ Z\Sigma V^H$ and $X_1$ (the first $3$ columns of $X$) is set to $UZ^{-1}$, where $U$ and $V$ are randomly generated matrices with orthonormal columns, $\Sigma=\mbox{diag}(\sigma_1,\sigma_2,\sigma_3)= (1^2,0.5^2,0.01^2)$ and  $Z$ is an arbitrary nonsingular matrix\footnote{Note that the choice of $Z$ does not affect the condition numbers $\kappa_i(A,B,1)$.}. Elements of the last row of $Y$ and the last column of $X$ were taken as random numbers. 
Since we have  $X_1Y_1^H=U\Sigma V^H$, we have  $\kappa_1(A,B,1)=2,\kappa_2(A,B,1)=1$ and $\kappa_3(A,B,1)=0.058$, the same condition numbers as the above second case with non-Hermitian perturbation, as was intended. The perturbations $E$ and $F$ are taken as arbitrary non-Hermitian matrices. 

In summary we tested under three different situations, all of which have the same $\sigma_i(X_1Y_1^H)$: (i) Both $(A,B)$ and $(E,F)$ are Hermitian (shown as ``Her + Her'' in Table \ref{exresults}), (ii)  $(A,B)$ is Hermitian  but $(E,F)$ is non-Hermitian (``Her + NonHer''), and (iii) Both $(A,B)$ and $(E,F)$ are non-Hermitian (``NonHer + NonHer''). 

The results are summarized in Table \ref{exresults} below, which shows the average  and maximum (shown as avg. and max respectively) values of
 $\Delta\lambda_i/\epsilon=|\widehat \lambda_i-1|/\epsilon$ among the  $10^6$ runs with randomly generated $E$ and $F$, along with the condition numbers  $\kappa_i(A,B,\lambda_0)$ (which are first order upper bounds for $\Delta\lambda_i/\epsilon$) for $i=1,2,3$.
\begin{table}[htbp]
\begin{center}
\caption{Average and maximum perturbation $\Delta\lambda_i/\epsilon$ of $10^6$ runs for $i=1,2,3$.}
\label{exresults}
    \begin{tabular}{|c|c|c|c||c|c|c||c|c|c|}\hline
& \multicolumn{3}{|c||}{Her + Her}& \multicolumn{3}{|c||}{Her + NonHer}& \multicolumn{3}{|c|}{NonHer + NonHer}\\\hline
$i$& avg.&    $\max$&$\kappa_i$ & avg. &    $\max$ &$\kappa_i$ & avg. &    $\max$&$\kappa_i$\\\hline
$1$&0.579& 1.98& 2.0   &0.41& 1.86 & 2.0   &0.42&1.90& 2.0  \\
$2$&0.141& 0.84& 1.0  &0.136& 0.76 & 1.0   &0.137&0.76& 1.0  \\
$3$&0.00018& 0.012& 0.02&0.00021& 0.027& 0.058 &0.00021&0.027 & 0.058\\\hline
    \end{tabular}
\end{center}
  \end{table}

We make the following observations. 
\begin{itemize}
\item We confirm that $\kappa_i$ is an upper bound of $\max\Delta\lambda_i/\epsilon$ for all $i$ in all three cases (which is necessarily true in the limit $\epsilon\rightarrow 0$). Interestingly, for $i=1$ the bound $\kappa_i$ is nearly attained while for $i=2,3$, $\max\Delta\lambda_i/\epsilon$ is noticeably smaller than $\kappa_i$, which suggests that for larger $i$ it becomes more and more rare that the largest-possible perturbation is attained. 
\item Reflecting the fact that $\kappa_i$ are the same for all the three cases for $i=1$ and $2$, 
we can see that $\max\Delta\lambda_i/\epsilon$ are similar in all three cases, so two eigenvalues have similar maximum sensitivities regardless of whether $A,B,E,F$ are Hermitian or not. 
On the contrary, $\max\Delta\lambda_i/\epsilon$ for $i=3$ show the somewhat different sensitivities of the third eigenvalue  depending on the structure of $E,F$. 
\item The behavior of the multiple eigenvalue is remarkably similar for the latter two cases, not only in terms of $\max\Delta \lambda_i/\epsilon$ but also avg.$\Delta \lambda_i/\epsilon$. 
This reflects the fact that the first order expansions of $\lambda_0$ are the same for the two cases, so that 
 the local behavior of an eigenvalue is determined solely by the singular values of $X_1Y_1^H$, and does not depend on the structure of the matrices $A$ and $B$. 
 \item Comparing avg.$\Delta\lambda_i/\epsilon$ with $\max\Delta\lambda_i/\epsilon$, 
 we see that the former is much smaller than the latter for larger $i$. For $i=1$ the difference seems less significant.
\end{itemize}
A precise explanation for the last two observations, 
which necessarily involves statistical analysis, 
 is an open problem: our discussions deal only with the maximum attainable perturbation $\max\Delta \lambda_i/\epsilon$, not with avg.$\Delta \lambda_i/\epsilon$. 
\end{example}
\subsection{Defective and infinite cases}\label{defsec}
So far we have treated only the case where $\lambda_0$ is a finite and nondefective multiple eigenvalue. Here we briefly consider the cases where $\lambda_0$ is infinite and/or defective.

The case $\lambda_0=\infty$ can be treated as in \cite{1storder,holder09} simply by considering the multiple zero eigenvalue of the pair $(B,A)$, for which the exact same discussion as above is valid.

When $\lambda_0$ is defective,  Lidskii's perturbation theory \cite{1storder,holder09} shows that the leading term in $\lambda_0$'s perturbation expansion  is not linear in $\epsilon$. Specifically, if $\lambda_0$ is an eigenvalue of $(A,B)$ of multiplicity $n_1r$ belonging to a Jordan block of dimension $n_1$ repeated $r$ times, then there are $n_1r$ eigenvalues of $(A+\epsilon E,B+\epsilon F)$ admitting the expansion
\begin{equation}
  \label{defcase}
\widehat \lambda_{i,\ell}= \lambda_0 +\left(\lambda_i(Y_1^H(E-\lambda_0F)X_1)\right)^{1/n_1}\epsilon^{1/n_1}+o(\epsilon^{1/n_1})
\end{equation}
for $i=1,2,\ldots, r$ and $\ell=1,2,\ldots, n_1$. 
Here $Y_1^H\in\mathbb{C}^{r\times n}$ and $X_1\in\mathbb{C}^{n\times r}$ represent the linearly independent left and right eigenvectors of $(A,B)$ corresponding to $\lambda_0$, and the value $\left(\lambda_i(Y_1^H(E-\lambda_0F)X_1)\right)^{1/n_1}$ takes all the $n_1$ distinct $n_1$th roots. 

We observe in \eqref{defcase} that although the leading exponent is different from that in \eqref{expand}, 
 the sensitivities of the multiple eigenvalue are still governed by  $|\lambda_i(Y_1^H(E-\lambda_0F)X_1)|$ for $i=1,\ldots,r$, for which we gave a bound  in the above discussions. Hence all our previous results carry over to this case, and 
the condition numbers of  $\lambda_0$  with the exponent $1/n_1$,  which we define by the theoretical bounds for
 $\sup_{\|E\|_2\leq 1,\|F\|_2\leq \tau}|\lambda_i(Y_1^H(E-\lambda_0F)X_1)|^{1/n_1}$,
are
\[\kappa_{i,n_1}(A,B,\lambda_0)=\left((1+\tau|\lambda_0|)\left(\prod_{j=1}^{i}\sigma_j\right)^{1/i}\right)^{1/n_1},\quad i=1,\ldots,r.\]
 By \eqref{defcase}, we must have $|\widehat \lambda_{i,\ell}-\lambda_0|/\epsilon^{1/n_1}\leq \kappa_{i,n_1}(A,B,\lambda_0)$ for $i=1,\ldots,r$ in the limit $\epsilon\rightarrow 0$ for any $E$ and $F$. Note that $\kappa_{i,n_1}(A,B,\lambda_0)$ does not depend explicitly on $\ell$. 
We observe that in the defective case $n_1\geq 2$, the exponent $1/n_1$ makes 
the difference among the condition numbers less significant than in the nondefective case. See the example below for an illustration. 
\subsubsection{Example}\label{ex2}
To examine the behavior of a defective multiple eigenvalue, we generate a $7$-by-$7$ pair $(A,B)$ defined by 
\begin{equation}  \label{eq:defective}
A=Y^{-H}\begin{bmatrix}
J& & &\\
&J & &\\
&  &J& \\
&  && 2
\end{bmatrix}X^{-1},\quad\mbox{and}\quad
B=Y^{-H}X^{-1},
\end{equation}
where $J=\begin{bmatrix}1&1\\0 &1 \end{bmatrix}$ is a $2$-by-$2$ Jordan block.
 $(A,B)$ has a multiple eigenvalue $\lambda_0=1$ of multiplicity six and a simple eigenvalue $2$. 
 $Y_1^H\equiv [Y(:,2)\ Y(:,4)\ Y(:,6)]^H=Z\Sigma V^H$ and 
$X_1\equiv [X(:,1)\ X(:,3)\ X(:,5)]=UZ^{-1}$ are the left and right eigenvectors corresponding to $\lambda_0$, where $U$ and $V$ are random matrices with orthonormal columns and $Z$ is an arbitrary nonsingular matrix. The other rows of $Y^H$ and columns of $X$ 
do not affect the condition numbers of $\lambda_0$, so we let them take random values. We let $\Sigma=\mbox{diag}(1^2,0.5^2,0.01^2)$, so that $\sigma_i(X_1Y_1^H)$ take the same values as in the non-Hermitian case of the second example in section \ref{ex1}.

Recall from \eqref{defcase} that perturbation in $(A,B)$  generally makes $\lambda_0$ split into $n_1r$ perturbed eigenvalues $\widehat\lambda_{i,\ell}$ for $i=1,\ldots,r$ and $\ell=1,\ldots,n_1$. 
 \eqref{defcase} also shows that for a fixed $i$, $|\widehat\lambda_{i,\ell}-\lambda_0|$ must be nearly equal for all $\ell$  up to $o(\epsilon^{1/n_1})$. 
For the matrix pair \eqref{eq:defective} we have $r=3$  and $n_1=2$, so we separate the six eigenvalues $\widehat\lambda_{i,\ell}$ into three groups according to the value of $i$, so that the two eigenvalues of the $i$th group have perturbation sensitivity governed by $\left|\lambda_i(Y_1^H(E-\lambda_0F)X_1)\right|^{1/n_1}$. 

With $\tau=1$, the condition numbers $\kappa_{i,2}(A,B,1)$ for the $i$th group for $i=1,2,3$ are
 $\kappa_{1,2}(A,B,1)=(2\cdot 1)^{1/2}=\sqrt{2},\kappa_{2,2}(A,B,1)=\left(2\cdot(1\cdot 0.5^2)^{1/2}\right)^{1/2}=1$ 
and $\kappa_{3,2}(A,B,1)=\left(2\cdot(1\cdot 0.5^2\cdot 0.01^2)^{1/3}\right)^{1/2}\simeq 0.24$. 
Comparing these with $\kappa_{i}(A,B,1)$ in the example in section \ref{ex1} we see that that although $\sigma_i(X_1Y_1^H)$ take the same values, the relative difference among the condition numbers is smaller here, due to the exponent $1/2$. 

Recalling that we must have $|\widehat \lambda_{i,\ell}-1|/\epsilon^{1/2}\leq \kappa_{i,2}(A,B,1)$ for small $\epsilon$, 
here we examine how large $|\widehat\lambda_{i,\ell}-1|/\epsilon^{1/2}$ becomes for $i=1,2,3$. To do this, of the six eigenvalues of $(A+\epsilon E,B+\epsilon F)$ close to $\lambda_0$, we check the perturbation of the most perturbed, third perturbed, and the fifth perturbed ones. 

In the experiment we let $E,F$ be random non-Hermitian matrices with $||E||_2,||F||_2\leq 1$, let $\epsilon=10^{-6}$ and tested with $10^6$  pairs. In table \ref{exresults2} we report the average  and maximum values of  $|\widehat \lambda_{i,\ell}-1|/\epsilon^{1/2}$ for $i=1,2,3$. 
\begin{table}[htbp]
\begin{center}
\caption{Defective matrix pair \eqref{eq:defective} with three $2\times 2$ Jordan blocks, average and maximum perturbation $|\widehat \lambda_{i,\ell}-1|/\epsilon^{1/2}$ of $10^6$ runs for $i=1,2,3$.}
\label{exresults2}
    \begin{tabular}{|c|c|c|c|}\hline
$i$& avg.&max & $\kappa_{i,2}$\\\hline
$1$&    0.511 &   1.21   &1.41\\
$2$&    0.282 &   0.733   &1\\
$3$&    0.0089 &   0.138   &0.24\\\hline
    \end{tabular}
\end{center}
  \end{table}

Similarly to the experiments in section \ref{ex1} for nondefective multiple eigenvalues, we see that a defective multiple eigenvalue also exhibits different sensitivities under perturbation. 
We also tested with Hermitian perturbations $E=E^H$ and $F=F^H$, and obtained nearly 
 the same results as in Table \ref{exresults}. This suggests that the structure of the perturbation does not play a role here. 

In all our experiments we had $|\widehat \lambda_{i,1}-\widehat \lambda_{i,2}|/\epsilon^{1/2}< 0.04$ for $i=1,2,3$, which matches the theoretical result indicated by \eqref{defcase} that for a given $i$, $|\widehat\lambda_{i,\ell}-1|$ are equal up to $o(\epsilon^{1/n_1})$ for all $\ell$. 

Finally, a comparison between Table \ref{exresults2} and the third case of Table \ref{exresults}  suggests that the relative difference among the multiple eigenvalues is smaller in the defective case, reflecting the last remark before this example. 
\begin{acknowledgements}
I thank the referees for their helpful suggestions, which improved the paper considerably. 
\end{acknowledgements}
\bibliographystyle{spmpsci}




\end{document}